\documentclass[graybox]{svmult}

\usepackage[protrusion=true,expansion,kerning,spacing]{microtype}
\usepackage{type1cm}        % activate if the above 3 fonts are
\usepackage{makeidx}         % allows index generation
\usepackage{graphicx}        % standard LaTeX graphics tool
\usepackage{multicol}        % used for the two-column index
\usepackage[bottom]{footmisc}% places footnotes at page bottom

\usepackage{newtxtext}       % 
\usepackage[varvw]{newtxmath}       % selects Times Roman as basic font
\usepackage{listings}
\definecolor{OliveGreen}{HTML}{3C8031}
\usepackage{tikz,pgfplots}
\usepackage{tikz-network}
\pgfplotsset{compat=1.18}
\usepackage{nicefrac}
\usepackage[inline]{enumitem}

\makeindex             
\begin{document}

\title*{PSCToolkit: solving sparse linear systems with a large number of GPUs}
\author{Pasqua D'Ambra\orcidID{0000-0003-2047-4986} and\\ Fabio Durastante\orcidID{0000-0002-1412-8289} and\\ Salvatore Filippone\orcidID{0000-0002-5859-7538}}
\institute{Pasqua D'Ambra\at Consiglio Nazionale delle Ricerche, Istituto per le Applicazioni del Calcolo ``Mauro Picone'', Via Pietro Castellino 111, 80131 - Napoli, \email{paqsua.dambra@cnr.it}
\and Fabio Durastante \at Università di Pisa, Dipartimento di Matematica, Largo Bruno Pontecorvo, 5 56127 - Pisa \email{fabio.durastante@unipi.it\\Consiglio Nazionale delle Ricerche, Istituto per le Applicazioni del Calcolo ``Mauro Picone'', Via Pietro Castellino 111, 80131 - Napoli.}
\and Salvatore Filippone \at Università di Roma ``Tor Vergata'', Dipartimento di Ingegneria Civile e Ingegneria informatica, Via del Politecnico, 1 00133 - Roma, \email{salvatore.filippone@uniroma2.it}\\Consiglio Nazionale delle Ricerche, Istituto per le Applicazioni del Calcolo ``Mauro Picone'', Via Pietro Castellino 111, 80131 - Napoli.
}
\authorrunning{P. D'Ambra, F. Durastante and S. Filippone}

\maketitle

\abstract*{In this chapter, we describe the Parallel Sparse Computation Toolkit (PSCToolkit), a suite of libraries for solving large-scale linear algebra problems in an HPC environment. In particular, we focus on the tools provided for the solution of symmetric and positive-definite linear systems using up to 8192 GPUs on the EuroHPC-JU Leonardo supercomputer. PSCToolkit is an ongoing mathematical software project aimed at exploiting the extreme computational speed of current supercomputers for relevant problems in Computational and Data Science. The toolkit is designed for node-level efficiency, flexibility and usability, supporting integration with both Fortran and C/C++, enabling researchers and developers from diverse computational backgrounds to leverage its powerful capabilities.
}

\abstract{In this chapter, we describe the Parallel Sparse Computation Toolkit (PSCToolkit), a suite of libraries for solving large-scale linear algebra problems in an HPC environment. In particular, we focus on the tools provided for the solution of symmetric and positive-definite linear systems using up to 8192 GPUs on the EuroHPC-JU Leonardo supercomputer. PSCToolkit is an ongoing mathematical software project aimed at exploiting the extreme computational speed of current supercomputers for relevant problems in Computational and Data Science. The toolkit is designed for flexibility and usability, supporting integration with both Fortran and C/C++, enabling researchers and developers from diverse computational backgrounds to leverage its powerful capabilities.
}

\section{Introduction}
\label{sec:introduction}

Solving large-scale linear systems is a fundamental challenge in many scientific and engineering disciplines, including
computational fluid dynamics~\cite{MR3425300,PSCToolkitAlya}, geosciences~\cite{PSCToolkitRichards}, and large-scale data
analysis~\cite{MR4072455}. These systems often involve millions or even billions of variables, requiring high computational power and
efficient algorithms for timely and accurate solutions. 
High-Performance Computing (HPC) environments provide the necessary
infrastructure, with Graphics Processing Units (GPUs) from vendors like NVIDIA, AMD, and Intel emerging as powerful
components that offer significant performance improvements over traditional CPU-based (Central--Processing--Unit-based) approaches~\cite{evodongarra}.

GPUs have a  massively parallel architecture which can execute thousands of threads concurrently, making them particularly
effective for parallelizing many linear algebra operations involving dense matrices, such as the BLAS-3 
\verb|GEMM|. Sparse matrix computations are
more challenging on GPUs due to their irregular memory accesses and lower arithmetic intensity. Therefore, innovative
algorithms and suitable data structures as well as implementation design patterns which leverage
memory locality and SIMD computations are needed to have substantial performance gains in solving sparse linear 
systems on GPUs.

This paper explores the methodologies implemented in the Parallel Sparse Computation Toolkit
(\texttt{PSCToolkit})~\cite{DAMBRA2023100463}, emphasizing the performance benefits of using NVIDIA GPUs for solving 
large-scale linear systems in HPC environments. We examine various algorithmic strategies developed specifically 
for GPU acceleration and present case studies and performance benchmarks to provide a comprehensive overview of the
state-of-the-art in GPU-accelerated linear system solvers. 
Our investigation demonstrates how \texttt{PSCToolkit} and GPUs
can significantly reduce computation time and enhance scalability, pushing the boundaries of what is achievable in
scientific and engineering computations.

In Sect.\ref{sec:psctoolkit}, we describe the \texttt{PSCToolkit} suite of libraries and their GPU-compatible
functionalities. Sect.\ref{sec:numerics} showcases performance results achieved on the EuroHPC-JU Leonardo supercomputer using up to
8192 GPUs. Finally, in Sect.~\ref{sec:conclusions}, we discuss the results and outline future developments for the suite,
both algorithmically and in implementation.

\vspace{-0.5em}\section{PSCToolkit: Parallel Sparse Computation Toolkit}
\label{sec:psctoolkit}

The Parallel Sparse Computation Toolkit (\texttt{PSCToolkit}) comprises three main components:
\begin{description}
\item[PSBLAS] \emph{Parallel Sparse BLAS} facilitates the parallelization of computationally intensive
scientific applications, focusing on the parallel implementation of iterative solvers for sparse
linear systems using the distributed-memory paradigm with message 
passing~\cite{FilipponeButtari,psblas2000}. It includes routines for sparse
matrix-dense matrix multiplication, block diagonal system solutions, and sparse matrix preprocessing,
 along with additional routines for dense matrix operations. Recent versions support large index
 spaces with 8-byte integers for extreme-scale problems.
\item[AMG4PSBLAS] \emph{Algebraic Multigrid Preconditioners for PSBLAS} provides parallel algebraic preconditioners,
evolving from the \texttt{MLD2P4} project~\cite{DambraDiSerafinoFilippone}.
Initially implementing additive-Schwarz type domain decomposition preconditioners and smoothed aggregation multilevel
hierarchies~\cite{BrezinaVanek,VanekMandelBrezina}, it now includes new smoothers and multigrid cycles, leveraging parallel
Coarsening based on Compatible Weighted Matching~\cite{MR4331965,DAmbraFilipponeVassilevski}.
\item[PSBLAS-EXT] \emph{PSBLAS extensions} offers additional matrix storage formats~\cite{Cardellini2014,FilipponeButtari},
and interfaces to external libraries \texttt{SPGPU} for NVIDIA GPU computations and \texttt{LIBRSB}~\cite{DBLP:conf/cata/MartoneFTPG10} for multi-core parallel machines.
\end{description}

The \texttt{PSCToolkit} is implemented in Fortran 2008~\cite{metcalf2008fortran}, reusing and adapting various existing
Fortran codes, along with some C and C++ routines. The library can be accessed from both Fortran and C/C++ at the
user level.

\subsection{PSBLAS applications}
\label{subsec:PSBLAS-applications}
\begin{figure}[htbp]
\sidecaption
\includegraphics[width=0.58\columnwidth]{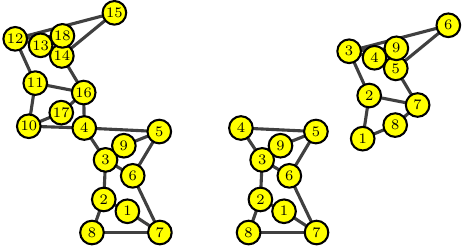}
\caption{Assignment of indexes in global and local numbering for two processes. The left panel represents the adjacency graph of the matrix considered with its global numbering, the two lateral graphs are the subgraphs relating to the locally numbered indices of the two processes.}
\label{fig:global_local}
\end{figure}

A generic application using \texttt{PSBLAS} involves several key steps. First, the parallel and accelerator
environments
are set up with \lstinline{psb_init} and \lstinline{psb_cuda_init}, followed by index space initialization with
\lstinline{psb_cdall}. Sparse matrices and dense vectors are allocated using \lstinline{psb_spall} and
\lstinline{psb_geall}, respectively. PSBLAS supports user-defined global to local index mappings, partitioning the index
space among processes and translating global numbering to local numbering for each process. The user normally defines
matrix entries  in the ``global'' numbering, and the library will translate them into the appropriate ``local'' one; see
Fig.~\ref{fig:global_local}.
Matrix and vector entries are generated and inserted with \lstinline{psb_spins} and \lstinline{psb_geins}, typically by
iterating over a finite element, finite volume mesh, or a finite difference stencil. These entities are then assembled
using \lstinline{psb_cdasb}, \lstinline{psb_spasb}, and \lstinline{psb_geasb}. PSBLAS uses an object-oriented model,
allowing users to handle sparse matrices and dense vectors independently of their ultimate storage format, which can be
chosen dynamically to suit the architecture without changing the application structure. 
This is achieved using so-called  \verb|mold| variables, following the techniques developed 
in~\cite{Cardellini2014};  see Fig.~\ref{fig:mold_code}.
\begin{figure}[t]
\begin{lstlisting}
type(psb_d_cuda_hlg_sparse_mat), target :: ahlg
type(psb_d_vect_cuda), target           :: dvgpu
type(psb_i_vect_cuda), target           :: ivgpu
call psb_cdasb(desc,info,mold=ivgpu)
call psb_spasb(a,desc,info,mold=ahlg)
call psb_geasb(b,desc,info,mold=dvgpu)
\end{lstlisting}
\caption{Example of usage of the \emph{mold} option to have data residing on the GPU. Specifically we require the GPU version of the  Hacked ELLPack format~\cite{Filippone:2017} with the corresponding dense GPU vector formats for the right-hand side \lstinline{b}, and the communicator index space \lstinline{desc}.}
\label{fig:mold_code}
\end{figure}
We stress that the library supports changing the storage format of a matrix object at runtime, allowing
efficient adaptation to different operation sequences with minimal user effort. After assembling the system, various
Krylov-type iterative methods~\cite{MR1990645} (e.g., Conjugate Gradient, GMRES, BiCGstab) can be used to 
approximate the solution of the linear system. These methods are well-suited for parallel environments, 
as they primarily require matrix-vector products which typically involve neighbourhood data communications 
among the various parallel processes. Note that, given that their arithmetic intensity is of order $\mathcal{O}(1)$, these algorithms  are memory/communication-bound algorithms; hence, as demonstrated by the HPCG benchmark results\footnote{\url{https://top500.org/lists/hpcg/2024/11/}}, their performance is limited by memory/communication bandwidth rather than by the floating-point throughput of processors. 

For an effective solution procedure we need a preconditioner, which transforms the linear system 
\begin{equation}\label{eq:linearsystem}
    A \mathbf{x} = \mathbf{b}, \qquad A \in \mathbb{K}^{N \times N}, \mathbf{x}, \mathbf{b} \in \mathbb{K}^{N}, \quad \mathbb{K} = \mathbb{R}, \mathbb{C},
\end{equation}
into one of the equivalent formulations
\[
B^{-1} A \mathbf{x} = B^{-1} \mathbf{b}, \quad \begin{cases}
A B^{-1} \mathbf{y} = \mathbf{b},\\
\mathbf{x} = B \mathbf{y},
\end{cases}
\begin{cases}
B_1^{-1} A B_2^{-1} \mathbf{y} = B_1^{-1} \mathbf{b},\\
\mathbf{x} = B_2 \mathbf{y},
\end{cases}
\]
with $B$, $B_1$, and $B_2$ being suitable invertible matrices that approximate $A$ and are simple to apply in a distributed environment using GPU accelerators. Advanced preconditioner tools are available in AMG4PSBLAS component.

\subsection{GPU support}
\label{subsec:PSBLAS-GPU}
As we mentioned, the storage format for the various sparse matrix objects is dynamic and can be 
chosen at  runtime using the techniques in~\cite{Cardellini2014}; in particular, this means that by
declaring some \verb|mold| variables and linking in the relevant code, it is possible to
extend the computational capabilities without touching the Krylov solver code, or indeed the structure
of the application. 

The use of these techniques allows us to introduce support for GPU computations through the use of appropriate \verb|mold| variables. The basic ideas are as follows:
\begin{itemize*}
    \item Every GPU-enabled data structure has a dual memory footprint, on the host and on the device side;
    \item For all data items, we keep track of whether it is the device or the host side that holds
    the latest   data updates;
    \item When data items in a method (such as matrix-vector product or AXPY) are GPU-enabled,
    computations are performed on the GPU by a kernel written in CUDA, and the output variables are
    marked as being updated on the device side;
    \item We synchronize the memory contents between the accelerator and the host only when absolutely necessary, since this involves a nontrivial  overhead.  
\end{itemize*}
The above rules enable the code implementing the Krylov solver to be executed seamlessly on the GPU
side. There are multiple possibilities for the sparse matrix storage formats, and the optimal choice
on the GPU side can be quite different from that on the CPU; for further details
see~\cite{Filippone:2017}. The GPU side code also includes some methods that handle packing and
unpacking of data for communication through MPI; this is kept to the barest minimum. 
One remaining issue with this architecture is the process of building the data structure: this is
still executed mostly on the CPU side, although some investigations are ongoing to explore
alternatives for improvement; this means that the best use scenario for the library is when the
matrices and/or preconditioners can be reused across multiple linear system solutions. To alleviate
this issue, we define interfaces that allow to update data coefficients in an existing matrix if the
structure is preserved; by the same token, we have split the Algebraic MultiGrid (AMG) preconditioner build phase in 
two steps, the construction of the matrix hierarchy, and the construction of the smoothers, 
the latter of which can be executed multiple times reusing an already assembled hierarchy.

Our toolkit supports many different storage formats on either CPU or GPU; in this paper we report
experimental data obtained with the use of the Hacked ELLPack format~\cite{Filippone:2017}
since we have found it provides  good performance in a wide variety of usage configurations. It is a
variation on the ELLPACK format; ELLPACK was designed to allow for vectorization, and its features 
blend well with the GPU architecture. The use of 2D arrays allows for regular access patterns, thus
ensuring good use of the GPU memory bandwidth; each row of the matrix is assigned to either one or two
threads, depending on the number of nonzeros, and the various threads in a warp access data in
coalesced memory transactions. 

The ELLPACK format has a drawback in that the use of a 2D array may require padding with zeros, and
thus it can  generate an excessive memory usage if there are some rows in the matrix with many more
nonzeros than the average. 
In our Hacked ELLPACK we partition the matrix into groups of rows, based on the size of a warp, and
store each group of rows as an independent ELLPACK matrix; in this way, any long row in the sparse
matrix will generate padding limited to a small subset of the rows, and the memory requirements do not
grow excessively.

\subsection{AMG4PSBLAS}
\label{subsec:amg4psblas}

AMG4PSBLAS is a package of parallel algebraic multilevel preconditioners. 
A general AMG method applied to a system of the form~\eqref{eq:linearsystem} is a stationary iterative 
method of the form
\[
\mathbf{x}^{(k)}=\mathbf{x}^{(k-1)}+B(\mathbf{b}-A\textbf{x}^{(k-1)}), \ \ k=1, 2, \ldots \ \ \text{given} \ \mathbf{x}^{(0)} \in \mathbb{K}^N,
\]
with matrix $B$ defined recursively through e.g. a $V-$cycle. More specifically,  we denote by $A_l$ the sequence of coarsened matrices 
computed via the triple-matrix Galerkin product:
\[
A_{l+1}=(P_l)^TA_lP_l, \ \ l=0, \ldots, \ell-1,
\]
with $A_0=A$ and $P_l$ a sequence of prolongation matrices of progressively smaller sizes $n_l \times n_{l+1}$, $n_{l+1} <
n_l$ and $n_0=n$, and by $M_l$ any $A_l$-convergent smoother for $A_l$, i.e., $\|I-M_l^{-1}A_l\|_{A_l} < 1$, where $I$
is the identity matrix of size $n_l$ and $\| \cdot \|_{A_l}$ indicates the $A_l$ norm. Then the matrix $M$ for the $V-$cycle, with the application of one sweep of pre- and post-smoothing, is the linear operator corresponding to the multiplicative composition of the following error propagation matrices:
\begin{equation}
\label{eq:amg}
I-B_lA_l= (I-(M_l)^{-T}A_l)(I-P_lB_{l+1}(P_l)^TA_l)(I-M_l^{-1}A_l) \ \ \forall l< \ell, 
\end{equation}
assuming that $B_{\ell} \approx A_{\ell}^{-1}$ is an approximation of the inverse of the coarsest-level matrix~\cite{V2008}. To fully specify the construction of the AMG preconditioner it is necessary to establish both the strategy with which to build the projectors and a reasonable choice for the smoothers. 
In the next two subsections we discuss some choices suitable for usage of the library on a large number of GPUs.

\subsubsection{Aggregation strategies}

The AMG4PSBLAS library provides two primary aggregation strategies with several variations:
\begin{description}
    \item[Decoupled Van\v{e}k, Mandel, Brezina.] This method groups fine grid nodes into disjoint sets called aggregates, each representing a single node in the coarse grid. Aggregates capture the problem's algebraic structure by identifying sets of strongly coupled nodes based on the user choice for $\theta$:
    \[
     N_i^{c}(\theta) = \{ j \mid |a_{i,j}| \geq \theta \sqrt{a_{i,i}a_{j,j}} \}, \quad \theta \in (0,1), \quad l=0,\ldots,n_{c}-1.
     \]
    Any remaining nodes are added to the nearest aggregates. The prolongation matrix is defined, 
    starting from a sample vector \(\mathbf{w}\) of the near kernel of \(A\) as:
     \begin{equation}\label{eq:prolongator_vbm}
         (\hat{P})_{i,j} = \begin{cases}
             w_i, & \text{if } i \in C_j,\\
             0, & \text{otherwise.}
         \end{cases}
     \end{equation}
    \item[Coupled Aggregation Based on Compatible Weighted Matching.] This method is independent of any 
    measure of  the strength of connection among nodes and finds is rationale in the Compatible Relaxation principles (see ~\cite{DAmbraFilipponeVassilevski} for details).
    Starting with the weighted adjacency graph \(G\) of \(A\) and edge weight matrix \(C\), it builds an approximate maximum weight matching \(\mathcal{M}\)~\cite{doi:10.1177/1094342012452893}. For each edge \(e_{i \mapsto j} \in \mathcal{M}\), vectors \(\mathbf{w}_e\) are identified as the orthonormal projection of a sample vector \(\mathbf{w}\) of the near kernel of \(A\) on the aggregate defined by \(\{i, j\}\). The prolongation matrix is:
\begin{equation}\label{eq:prolongator}
\hat{P} = \left (
\begin{array}{cc}
\tilde{P} & 0\\
0 & W
\end{array}
\right ) \in \mathbb{R}^{N \times n_c}, 
\end{equation}
where \(\tilde{P} = \operatorname{blockdiag}(\mathbf{w}_{e_1}, \ldots, \mathbf{w}_{e_{n_p}} )\) for the $n_p$ matched edges, and \(W = \operatorname{diag}(\mathbf{w}_s/|\mathbf{w}_s|)\) for the $n_s$ unmatched vertices, with \(n_c = n_p + n_s\). The range of \(\hat{P}\) includes \(\mathbf{w}\); to improve the size reduction, it is possible to perform multiple ($k$) matching sweeps, thereby generating aggregates of size $2^k$.
\end{description}

\begin{figure}[htbp]
\sidecaption
\includegraphics[width=3.2cm]{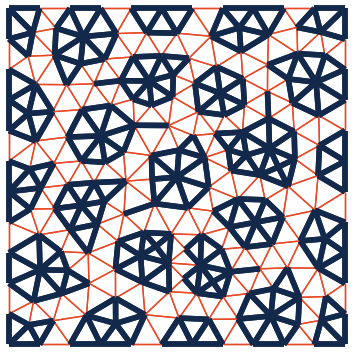}
\includegraphics[width=3.2cm]{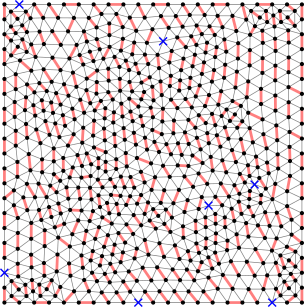}
\caption{Example of the aggregates obtained with the two different schemes. On the left, the aggregation of Van\v{e}k, Mandel and Brezina, on the right the aggregation obtained with the approximate matching of maximum weight.}
\label{fig:aggregates}       % Give a unique label
\end{figure}

In both strategies, the prolongation matrix is a piecewise constant interpolation operator. To enhance its regularity and robustness, a smoothing procedure is applied: \({P} = (I - \omega D^{-1}A)\hat{P}\), where \(D = \operatorname{diag}(A)\) and \(\omega = 1/\|D^{-1}A\|_{\infty} \approx 1/\rho(D^{-1}A)\), with \(\rho(D^{-1}A)\) being the spectral radius of \(D^{-1}A\). Fig.~\ref{fig:aggregates} illustrates examples of the two aggregation procedures on small matrices.
The preconditioner and its configuration can be implemented as simply as depicted in snippet of code in Fig.~\ref{fig:multigrid-definition}, in which we set an ML scheme with the coupled aggregation based on weighted graph matching.
\begin{figure}[htbp]
    \centering
\begin{lstlisting}
type(amg_dprec_type)  :: prec
call prec%init(ctxt,'ML',info) !Multigrid with
call prec%set('ml_cycle','VCYCLE',info) !V-Cycle
call prec%set('par_aggr_alg','COUPLED',info) !coupled
call prec%set('aggr_prol','SMOOTHED',info) !smoothed
call prec%set('aggr_type','MATCHBOXP',info) !matching-based
call prec%set('aggr_size',8,info) !aggregation (aggr. size le 8)
\end{lstlisting}
    \caption{Preconditioner instatiation and selection for a multigrid preconditioner using a V-cycle with the smoothed aggregation procedure based on the weighted graph matching procedure.}
    \label{fig:multigrid-definition}
\end{figure}
\vspace{-2em}\subsubsection{GPU enabled smoothers}

A GPU-friendly smoother primarily utilizes operations that can be efficiently implemented in a single
instruction multiple data (SIMD) architecture; these include smoothers that use matrix-vector products, \texttt{axpy}, and diagonal scaling operations~\cite{MR2861652,MR4331965,MR4663280}.
These smoothers include:
\begin{description}
    \item[$\ell_1$-Jacobi smoother.] This Jacobi method variant uses an additive correction to the diagonal. The
    smoother matrix at level $l$ is \[M_l = \operatorname{diag}(A_l) + \operatorname{diag}\left( \left\{
    \sum_{j=1, j \neq i}^{N} |a_{i,j}^{(l)}| \right\}_{i=1}^{N} \right).\] It guarantees convergence without
    needing a multiplicative scaling factor and is $A_l$-convergent~\cite{MR2861652}.
    \item[AINV smoother.] This method constructs an approximate inverse of \(A\) using a sparse matrix \(M_l
    \approx A_l^{-1}\)~\cite{MR1787297,MR3457372}. It involves an incomplete
    biconjugation, approximating \(A_l^{-1}\) as a product of two sparse triangular matrices, \(Z\)
    and \(W\), such that     \(A_l^{-1} \approx ZW^T\). In a distributed setting, it approximates the
    inversion of diagonal blocks in a   block-Jacobi smoother; on GPUs, it is efficient since it needs
    only to perform sparse matrix-vector products 
    with \(Z\) and \(W\).
    \item[INVK smoother.] This method approximates the inverse by inverting the sparse \(L\) and \(U\)
    factors of \(A\)~\cite{MR3457372,MR1699790}. Similar to AINV, it approximates the inversion of
    diagonal blocks in a block-Jacobi smoother, therefore it is also GPU-efficient, requiring only
    sparse matrix-vector products.
\end{description}
The selection for the  smoothers can be added to the configuration shown in 
Fig.~\ref{fig:multigrid-definition} with a few additional calls; after that, the preconditioner is built by
generating the aggregation hierarchy and the smoothers. 
\begin{figure}[h]
    \centering
\begin{lstlisting}
call prec%set('smoother_type','L1-JACOBI',info) !Smoother
call prec%set('smoother_sweeps',4,info) !4 iterations
call prec%set('coarse_solve','L1-JACOBI',info) !Coarse solver
call prec%set('coarse_sweeps', 30,info) !30 iterations
call prec%hierarchy_build(a,desc_a,info) ! build phase
call prec%smoothers_build(a,desc_a,info,amold=ahlg,vmold=dvgpu,imold=ivgpu)
\end{lstlisting}
\caption{Set up of the $\ell_1$-Jacobi method both as smoother (4 iterations) and as coarse solver (30 iterations), and calls to construct the hierarchy
and the smoother}% on the GPU side, through the use of the \texttt{'mold'} option.}
\label{fig:prec_build}
\end{figure}
\vspace{-1em}\section{Numerical Example}
\label{sec:numerics}

The experiments are performed on the CINECA Leonardo machine, ranked 9\textsuperscript{th} on the November 2024 TOP500 list\footnote{See~\url{https://top500.org/lists/top500/2024/11/}.}. Each node is a personalized Bull Sequana blade X2135 "Da Vinci," consisting of a single Intel Xeon 8358 CPU with 32 cores at 2.6 GHz, 512 GB of DDR4 RAM at 3200 MHz, and 4 NVIDIA custom Ampere A100 GPUs with 64 GB of HBM2 memory. All nodes are interconnected through an NVIDIA Mellanox network with DragonFly+, capable of a maximum bandwidth of 200 Gbit/s between each pair of nodes. The code is compiled with the \texttt{gcc/11.3.0} suite, using \texttt{openmpi/4.1.4} and CUDA \texttt{compilation tools, release 11.8, V11.8.89}. Comparisons are run with the NVIDIA AMGX library \texttt{v.2.4.0}~\cite{MR3504566}.

We test the solution of the 7-point finite difference discretization of the 3D Poisson problem:
\[
\begin{cases}
    -\nabla^2 u = 1, & \mathbf{x} \in [0,1]^3, \\
    u(\mathbf{x}) = 0, & \mathbf{x} \in \partial [0,1]^3,
\end{cases}
\]
using a synchronization-reduced variant of the Conjugate Gradient method~\cite{MR1797890}, aiming for a relative tolerance on the residual of $\tau = 10^{-6}$. To evaluate the weak-scaling properties of different algorithms, we use $8 \times 10^6$ unknowns per GPU, \emph{i.e.}, $3.2 \times 10^7$ unknowns per node. The largest system we consider has $ \approx 6.5 \times 10^{10}$ degrees of freedom. The call to the Krylov solver is again a single line to which we pass all the built quantities as depicted in Fig.~\ref{fig:fcg_call}.
\begin{figure}[htbp]
    \centering
\begin{lstlisting}
call psb_krylov('FCG',a,prec,b,x,1d-7,desc,info)
\end{lstlisting}
    \caption{Simplest possible call to the flexible conjugate gradient method; \lstinline{x} contains  the starting guess, and the approximate solution upon completion. Other control arguments are optional.}
    %\caption{Simplest possible call to the flexible conjugate gradient method. The vector \lstinline{x}  contains the starting guess at the beginning, and the approximate solution at the end. Other options like the maximum number of iteration, or how often to print residual estimate, are optional.}
     %can be passed to the same call.}
    \label{fig:fcg_call}
\end{figure}
\vspace{-2em}\subsection{Algorithmic tests}\label{subsec:tested-methods} %
After the discussion in Sect.~\ref{subsec:amg4psblas}, we tested the following preconditioners from  \texttt{PSCToolkit}:
\begin{description}
	\item[VBM]{decoupled Van\v{e}k, Mandel, Brezina aggregation~\cite{VanekMandelBrezina}, V-cycle, $\ell_1$-Jacobi smoother (4 sweeps), at most 40 iterations of the Preconditioned CG coupled to $\ell_1$-Jacobi preconditioner as coarsest solver;}
	\item[SMATCH]{matching-based aggregation with aggregates of maximum size equal to $8$, smoothing of prolongators, further algorithmic choices as in VBM;} 
	\item[VMATCH]{matching-based aggregation as in SMATCH, un-smoothed prolongators, Variable V-cycle\footnote{2 smoother iteration at the first level, and doubled at each following level.}, further algorithmic choices as in VBM.}
\end{description}
and we compare them against the following configurations from the NVIDIA AMGX library~\cite{MR3504566,AMGXusersguide}:
\begin{description}
    \item[AMGX CLASSICAL] coarsening done by classical, also known as Ruge-St\"{u}ben, AMG approach, where the coarse nodes are a subset of the fine nodes and a distance-2 interpolation is applied, V-cycle, smoother $\ell_1$-Jacobi (4 sweeps), $\ell_1$-Jacobi (40 sweeps) coarsest solver; 
    \item[AMGX AGGREGATION] aggregation by iterative parallel graph matching, aggregates of maximum size equal to $8$, further algorithmic choices as in AMGX CLASSICAL.
\end{description}
We begin by analyzing the algorithmic scaling performance of various multigrid methods. In
Fig.~\ref{fig:iteration_count}, we show the number of iterations required for the preconditioned
conjugate gradient to converge for each preconditioner.
\begin{figure}[htbp]
\centering
\definecolor{mycolor1}{rgb}{0.00000,0.44700,0.74100}%
\definecolor{mycolor2}{rgb}{0.85000,0.32500,0.09800}%
\definecolor{mycolor3}{rgb}{0.92900,0.69400,0.12500}%
\definecolor{mycolor4}{rgb}{0.49400,0.18400,0.55600}%
\definecolor{mycolor5}{rgb}{0.46600,0.67400,0.18800}%

\begin{tikzpicture}
\begin{axis}[%
	width=0.65\columnwidth,
	height=1.1in,
	at={(1.861in,0.689in)},
	scale only axis,
	xmode=log,
	xmin=1,
	xmax=8192,
	xtick={1,2,4,8,16,32,64,128,256,512,1024,2048,4096,8192},
	xticklabels={{1}, {2},{4},{8},{16},{32},{64},{128},{256},{512},{1024},{2048},{4096},{8192}},
	xminorticks=true,
	xlabel style={font=\color{white!15!black},at={(axis description cs:0.5,-0.1)}},
    ylabel style={font=\footnotesize},
    ylabel={Iterations},
	xticklabel style={rotate=35},
	ymode=log,
	ymin=9,
	ymax=282,
	yminorticks=true,
	axis background/.style={fill=white},
	title style={font=\bfseries},
	axis x line*=bottom,
	axis y line*=left,
	xmajorgrids,
	xminorgrids,
	ymajorgrids,
	yminorgrids,
	legend style={at={(1.0,0.5)}, anchor=west, legend cell align=left, align=left, draw=none, font=\scriptsize}
	]
	\addplot [color=mycolor1, dashed, line width=3.0pt, mark=o, mark options={solid, mycolor1}]
	table[row sep=crcr]{%
		1	18\\
		2	20\\
		4	21\\
		8	21\\
		16	23\\
		32	24\\
		64	24\\
		128	29\\
		256	32\\
		512	33\\
		1024	39\\
		2048	44\\
		4096	71\\
		8192	65\\
	};
	\addlegendentry{VBM}
	
	\addplot [color=mycolor2, dashed, line width=3.0pt, mark=x, mark options={solid, mycolor2}]
	table[row sep=crcr]{%
		1	9\\
		2	10\\
		4	11\\
		8	10\\
		16	11\\
		32	14\\
		64	10\\
		128	13\\
		256	18\\
		512	12\\
		1024	24\\
		2048	11\\
		4096	37\\
		8192	39\\
	};
	\addlegendentry{SMATCH}
	
	\addplot [color=mycolor3, dashed, line width=3.0pt, mark=triangle, mark options={solid, rotate=90, mycolor3}]
	table[row sep=crcr]{%
		1	28\\
		2	34\\
		4	37\\
		8	36\\
		16	48\\
		32	51\\
		64	49\\
		128	65\\
		256	69\\
		512	65\\
		1024	102\\
		2048	70\\
		4096	129\\
		8192	164\\
	};
	\addlegendentry{VMATCH}
	
	\addplot [color=mycolor4, line width=3.0pt, mark=star, mark options={solid, mycolor4}]
	table[row sep=crcr]{%
		1	62\\
		2	72\\
		4	82\\
		8	95\\
		16	110\\
		32	122\\
		64	131\\
		128	136\\
		256	146\\
		512	164\\
		1024	183\\
		2048	206\\
		4096	246\\
		8192	282\\
	};
	\addlegendentry{AMGX-AGGR.}
	
	\addplot [color=mycolor5, line width=3.0pt, mark=star, mark options={solid, mycolor5}]
	table[row sep=crcr]{%
		1	10\\
		2	10\\
		4	10\\
		8	10\\
		16	10\\
		32	11\\
		64	11\\
		128	11\\
		256	11\\
		512	12\\
		1024	12\\
		2048	12\\
		4096	15\\
		8192	18\\
	};
	\addlegendentry{AMGX-CLAS.}
	
\end{axis}
\end{tikzpicture}
    \caption{Total number of FCG iterations vs number of GPUs for the methods of Sec.~\ref{subsec:tested-methods}.}
    \label{fig:iteration_count}
\end{figure}
All versions show a moderate increase in the number of iterations; our proposed methods exhibit scalability
intermediate between the two preconditioners implemented in AMGX.
To better understand the iteration behaviour, we show in
Fig.~\ref{fig:op_complexity}, the operator complexity of all the AMG hierarchies built by the different methods. We can observe that the very good algorithmic scalability of the AMGX-CLASSICAL preconditioner is the result of a large operator complexity which indeed is larger than $4$, while, as expected, our VMATCH and AMGX-AGGREGATION show
the best complexities, not larger than $1.14$ and $1.33$, respectively. This is due to the sparsity of their prolongator operators.
\begin{figure}
    \centering
    \definecolor{mycolor1}{rgb}{0.00000,0.44700,0.74100}%
\definecolor{mycolor2}{rgb}{0.85000,0.32500,0.09800}%
\definecolor{mycolor3}{rgb}{0.92900,0.69400,0.12500}%
\definecolor{mycolor4}{rgb}{0.49400,0.18400,0.55600}%
\definecolor{mycolor5}{rgb}{0.46600,0.67400,0.18800}%
\definecolor{mycolor6}{rgb}{0.63500,0.07800,0.18400}%
    \begin{tikzpicture}
    
    \begin{axis}[%
width=0.7\columnwidth,
height=1.1in,
at={(1.861in,0.689in)},
scale only axis,
xmode=log,
xmin=1,
xmax=8192,
xtick={1,2,4,8,16,32,64,128,256,512,1024,2048,4096,8192},
xticklabels={{1}, {2},{4},{8},{16},{32},{64},{128},{256},{512},{1024},{2048},{4096},{8192}},
xminorticks=true,
xticklabel style={rotate=35},
xlabel style={font=\color{white!15!black},at={(axis description cs:0.5,-0.1)}},
ymin=1,
ymax=5,
yminorticks=true,
axis background/.style={fill=white},
ylabel style={font=\footnotesize},
ylabel={Operator Complexity},
axis x line*=bottom,
axis y line*=left,
xmajorgrids,
xminorgrids,
ymajorgrids,
yminorgrids,
legend style={at={(1.0,0.5)}, anchor=west, legend cell align=left, align=left, draw=none, font=\scriptsize}
]
\addplot [color=mycolor1, line width=3.0pt, mark=o, mark options={solid, mycolor1}]
  table[row sep=crcr]{%
1	1.575\\
2	1.578\\
4	1.58\\
8	1.583\\
16	1.584\\
32	1.584\\
64	1.587\\
128	1.588\\
256	1.587\\
512	1.589\\
1024	1.588\\
2048	1.59\\
4096	1.588\\
8192	1.588\\
};
\addlegendentry{VBM}

\addplot [color=mycolor2, dashed, line width=3.0pt, mark=x, mark options={solid, mycolor2}]
  table[row sep=crcr]{%
1	1.894\\
2	1.905\\
4	1.915\\
8	1.917\\
16	1.925\\
32	1.93\\
64	1.93\\
128	1.936\\
256	1.905\\
512	1.937\\
1024	1.942\\
2048	1.939\\
4096	1.906\\
8192	1.945\\
};
\addlegendentry{SMATCH} % VCYCLE-HLG

\addplot [color=mycolor3, dashed, line width=3.0pt, mark=triangle, mark options={solid, rotate=90, mycolor3}]
  table[row sep=crcr]{%
1	1.142\\
2	1.142\\
4	1.143\\
8	1.142\\
16	1.143\\
32	1.143\\
64	1.143\\
128	1.143\\
256	1.144\\
512	1.143\\
1024	1.144\\
2048	1.143\\
4096	1.144\\
8192	1.144\\
};
\addlegendentry{VMATCH}

\addplot [color=mycolor4, line width=3.0pt, mark=star, mark options={solid, mycolor4}]
  table[row sep=crcr]{%
1	1.27979\\
2	1.31187\\
4	1.33117\\
8	1.33162\\
16	1.32133\\
32	1.31887\\
64	1.31914\\
128	1.31421\\
256	1.31314\\
512	1.31329\\
1024	1.31091\\
2048	1.31041\\
4096	1.31049\\
8192	1.30932\\
};
\addlegendentry{AMGX-AGGR.}

\addplot [color=mycolor5, line width=3.0pt, mark=star, mark options={solid, mycolor5}]
  table[row sep=crcr]{%
1	4.45456\\
2	4.43576\\
4	4.51377\\
8	4.52376\\
16	4.51239\\
32	4.49595\\
64	4.50135\\
128	4.49925\\
256	4.49252\\
512	4.4952\\
1024	4.49503\\
2048	4.4921\\
4096	4.49354\\
8192	4.49371\\
};
\addlegendentry{AMGX-CLAS.}

\end{axis}
    \end{tikzpicture}
    
    \caption{Operator complexity obtained from the different coarsening schemes, that is $\nicefrac{\sum_{l=0}^{\ell}\operatorname{nnz}(A_l)}{\operatorname{nnz}(A_0)}$ which is both a measure of the memory occupation of the AMG preconditioner and a measure of its application cost. }
    \label{fig:op_complexity}
\end{figure}
To further investigate the efficiency of the
methods, we also analyze the solve time (in seconds) across different configurations, as shown in
Fig.~\ref{fig:solve-time}.
\begin{figure}[htbp]
    \centering
    \definecolor{mycolor1}{rgb}{0.00000,0.44700,0.74100}%
\definecolor{mycolor2}{rgb}{0.85000,0.32500,0.09800}%
\definecolor{mycolor3}{rgb}{0.92900,0.69400,0.12500}%
\definecolor{mycolor4}{rgb}{0.49400,0.18400,0.55600}%
\definecolor{mycolor5}{rgb}{0.46600,0.67400,0.18800}%

\begin{tikzpicture}
	\begin{axis}[%
	width=0.7\columnwidth,
	height=1.1in,
	at={(1.861in,0.689in)},
	scale only axis,
	xmode=log,
	xmin=1,
	xmax=8192,
	xtick={1,2,4,8,16,32,64,128,256,512,1024,2048,4096,8192},
	xticklabels={{1}, {2},{4},{8},{16},{32},{64},{128},{256},{512},{1024},{2048},{4096},{8192}},
	xminorticks=true,
	xlabel style={font=\color{white!15!black}},
    ylabel style={font=\footnotesize},
    ylabel={Solve time},
	xticklabel style={rotate=35},
	xlabel style={font=\color{white!15!black},at={(axis description cs:0.5,-0.1)}},
	ymode=log,
	ymin=0.14755,
	ymax=76.5363,
	yminorticks=true,
	axis background/.style={fill=white},
	axis x line*=bottom,
	axis y line*=left,
	xmajorgrids,
	xminorgrids,
	ymajorgrids,
	yminorgrids,
	legend style={at={(1.0,0.5)}, anchor=west, legend cell align=left, align=left, draw=none, font=\scriptsize}
	]
	\addplot [color=mycolor1, dashed, line width=3.0pt, mark=o, mark options={solid, mycolor1}]
	table[row sep=crcr]{%
		1	0.22308\\
		2	0.32834\\
		4	0.41095\\
		8	0.46954\\
		16	0.63889\\
		32	0.73148\\
		64	0.87132\\
		128	1.5738\\
		256	1.85849\\
		512	1.23709\\
		1024	2.49795\\
		2048	2.64294\\
		4096	5.19639\\
		8192	5.39526\\
	};
	\addlegendentry{VBM}
	
	\addplot [color=mycolor2, dashed, line width=3.0pt, mark=x, mark options={solid, mycolor2}]
	table[row sep=crcr]{%
		1	0.14755\\
		2	0.22388\\
		4	0.41069\\
		8	0.34686\\
		16	0.53371\\
		32	0.94287\\
		64	0.91508\\
		128	1.50342\\
		256	2.90102\\
		512	2.08193\\
		1024	4.77418\\
		2048	2.06253\\
		4096	9.90256\\
		8192	7.79931\\
	};
	\addlegendentry{SMATCH}
	
	\addplot [color=mycolor3, dashed, line width=3.0pt, mark=triangle, mark options={solid, rotate=90, mycolor3}]
	table[row sep=crcr]{%
		1	0.43282\\
		2	0.84559\\
		4	1.29211\\
		8	1.60107\\
		16	2.94846\\
		32	4.00064\\
		64	6.92191\\
		128	8.09396\\
		256	10.2308\\
		512	12.1772\\
		1024	22.0979\\
		2048	13.0986\\
		4096	76.5363\\
		8192	34.8665\\
	};
	\addlegendentry{VMATCH}
	
	\addplot [color=mycolor4, line width=3.0pt, mark=star, mark options={solid, mycolor4}]
	table[row sep=crcr]{%
		1	0.816997\\
		2	1.46825\\
		4	1.74968\\
		8	2.10821\\
		16	2.52544\\
		32	3.14764\\
		64	5.14636\\
		128	7.17638\\
		256	8.18452\\
		512	11.3246\\
		1024	37.8813\\
		2048	13.9117\\
		4096	13.8165\\
		8192	26.4272\\
	};
	\addlegendentry{AMGX-AGGR.}
	
	\addplot [color=mycolor5, line width=3.0pt, mark=star, mark options={solid, mycolor5}]
	table[row sep=crcr]{%
		1	0.300799\\
		2	0.385396\\
		4	0.405538\\
		8	0.467566\\
		16	0.511307\\
		32	0.728832\\
		64	1.30343\\
		128	1.9389\\
		256	3.4067\\
		512	4.86612\\
		1024	4.30711\\
		2048	5.39214\\
		4096	13.1165\\
		8192	12.1541\\
	};
	\addlegendentry{AMGX-CLAS.}
	
\end{axis}
\end{tikzpicture}%
    \caption{Solve time vs number of GPUs. Time needed for the FCG to reach convergence for the methods described in Sect.~\ref{subsec:tested-methods}.}
    \label{fig:solve-time}
\end{figure}
The results indicate that both VBM and SMATCH methods achieve a better solution time than the AMGX preconditioners, showing that, although the
small operator complexity produces a degradation effect on algorithmic scalability, the reduced application costs make our preconditioners more efficient for increasing number of GPUs, where the increase in the number of iterations is balanced by the smaller cost of each iteration (see Fig.~\ref{fig:time-per-iteration}).
\begin{figure}[htbp]
    \centering
    
    	\definecolor{mycolor1}{rgb}{0.00000,0.44700,0.74100}%
	\definecolor{mycolor2}{rgb}{0.85000,0.32500,0.09800}%
	\definecolor{mycolor3}{rgb}{0.92900,0.69400,0.12500}%
	\definecolor{mycolor4}{rgb}{0.49400,0.18400,0.55600}%
	\definecolor{mycolor5}{rgb}{0.46600,0.67400,0.18800}%
	
\begin{tikzpicture}
	\begin{axis}[%
	width=0.7\columnwidth,
height=1.1in,
at={(1.861in,0.689in)},
scale only axis,
xmode=log,
xmin=1,
xmax=8192,
xtick={1,2,4,8,16,32,64,128,256,512,1024,2048,4096,8192},
xticklabels={{1}, {2},{4},{8},{16},{32},{64},{128},{256},{512},{1024},{2048},{4096},{8192}},
xminorticks=true,
xlabel style={font=\color{white!15!black}},
ylabel style={font=\footnotesize},
ylabel={Time per iteration},
xticklabel style={rotate=35},
xlabel style={font=\color{white!15!black},at={(axis description cs:0.5,-0.1)}},
ymode=log,
ymin=0.01239,
ymax=1.2368,
yminorticks=true,
axis background/.style={fill=white},
axis x line*=bottom,
axis y line*=left,
xmajorgrids,
xminorgrids,
ymajorgrids,
yminorgrids,
legend style={at={(1.0,0.5)}, anchor=west, legend cell align=left, align=left, draw=none, font=\scriptsize}
		]
		\addplot [color=mycolor1, dashed, line width=3.0pt, mark=o, mark options={solid, mycolor1}]
		table[row sep=crcr]{%
			1	0.01239\\
			2	0.01642\\
			4	0.01957\\
			8	0.02236\\
			16	0.02778\\
			32	0.03048\\
			64	0.03631\\
			128	0.05427\\
			256	0.05808\\
			512	0.03749\\
			1024	0.06405\\
			2048	0.06007\\
			4096	0.07319\\
			8192	0.083\\
		};
		\addlegendentry{VBM}
		
		\addplot [color=mycolor2, dashed, line width=3.0pt, mark=x, mark options={solid, mycolor2}]
		table[row sep=crcr]{%
			1	0.01639\\
			2	0.02239\\
			4	0.03734\\
			8	0.03469\\
			16	0.04852\\
			32	0.06735\\
			64	0.09151\\
			128	0.11565\\
			256	0.16117\\
			512	0.1735\\
			1024	0.19892\\
			2048	0.1875\\
			4096	0.26764\\
			8192	0.19998\\
		};
		\addlegendentry{SMATCH}
		
		\addplot [color=mycolor3, dashed, line width=3.0pt, mark=triangle, mark options={solid, rotate=90, mycolor3}]
		table[row sep=crcr]{%
			1	0.01546\\
			2	0.02487\\
			4	0.03492\\
			8	0.04447\\
			16	0.06143\\
			32	0.07844\\
			64	0.14126\\
			128	0.12452\\
			256	0.14827\\
			512	0.18734\\
			1024	0.21665\\
			2048	0.18712\\
			4096	0.5933\\
			8192	0.2126\\
		};
		\addlegendentry{VMATCH}
		
		\addplot [color=mycolor4, line width=3.0pt, mark=star, mark options={solid, mycolor4}]
		table[row sep=crcr]{%
			1	0.0131774\\
			2	0.0203923\\
			4	0.0213376\\
			8	0.0221916\\
			16	0.0229586\\
			32	0.0258003\\
			64	0.0392852\\
			128	0.0527675\\
			256	0.0560584\\
			512	0.0690525\\
			1024	0.207002\\
			2048	0.0675326\\
			4096	0.0561645\\
			8192	0.0937135\\
		};
		\addlegendentry{AMGX-AGGR.}
		
		\addplot [color=mycolor5, line width=3.0pt, mark=star, mark options={solid, mycolor5}]
		table[row sep=crcr]{%
			1	0.0300799\\
			2	0.0385396\\
			4	0.0405538\\
			8	0.0467566\\
			16	0.0511307\\
			32	0.0662574\\
			64	0.118493\\
			128	0.176263\\
			256	0.3097\\
			512	0.40551\\
			1024	0.358925\\
			2048	0.449345\\
			4096	0.874433\\
			8192	0.675226\\
		};
		\addlegendentry{AMGX-CLAS.}

	\end{axis}
	
	\end{tikzpicture}
    
    \caption{Time per iteration vs number of GPUs. Average time needed to perform a single preconditioned iteration of the CG method with the various preconditioners, measured in seconds.}
    \label{fig:time-per-iteration}
\end{figure}

For the configuration exhibiting the best weak scalability performance in PSCToolkit, identified by the VBM label, we briefly examine its strong scalability properties. Specifically, we analyze the solution time by assigning a problem of size $480^3 \approx 1.1 \times 10^8$ to a single GPU. This corresponds to a matrix occupying 9.5 GB when stored in HLG format. The number of GPUs is then scaled from 1 to 16, resulting in a local problem size of approximately $6.9 \times 10^6$ unknowns for the 16-GPU case, comparable to the setup used for weak scaling analysis.
\begin{figure}[htbp]
    \centering

\definecolor{mycolor1}{rgb}{0.00000,0.44700,0.74100}%
\begin{tikzpicture}

\begin{axis}[%
width=0.328\columnwidth,
height=0.334\columnwidth,
at={(0\columnwidth,0\columnwidth)},
scale only axis,
xmode=log,
xmin=1,
xmax=16,
xtick={1,2,4,8,16},
xticklabels={{1},{2},{4},{8},{16}},
ytick={4,1,0.5},
yticklabels={4,1,0.5},
xminorticks=true,
ymode=log,
ymin=0.299815625,
ymax=4.79705,
yminorticks=true,
ylabel style={font=\color{white!15!black}},
ylabel={Solve time},
axis background/.style={fill=white}
]
\addplot [color=mycolor1, line width=2.0pt, mark=o, mark options={solid, mycolor1}, forget plot]
  table[row sep=crcr]{%
1	4.79705\\
2	2.42424\\
4	1.3777\\
8	0.904502\\
16	0.684951\\
};
\addplot [color=black, dashed, line width=2.0pt, forget plot]
  table[row sep=crcr]{%
1	4.79705\\
2	2.398525\\
4	1.1992625\\
8	0.59963125\\
16	0.299815625\\
};
\end{axis}

\begin{axis}[%
width=0.328\columnwidth,
height=0.334\columnwidth,
at={(0.48\columnwidth,0\columnwidth)},
scale only axis,
xmin=0,
xmax=16,
xtick={1,2,4,8,16},
ytick={1,10},
ymode=log,
ymin=1,
ymax=16,
yminorticks=true,
ylabel style={font=\color{white!15!black}},
ylabel={Speedup},
axis background/.style={fill=white}
]
\addplot [color=mycolor1, line width=2.0pt, mark=o, mark options={solid, mycolor1}, forget plot]
  table[row sep=crcr]{%
1	1\\
2	1.9787851037851\\
4	3.48192639907092\\
8	5.30352613924568\\
16	7.00349368056985\\
};
\addplot [color=black, dashed, line width=2.0pt, forget plot]
  table[row sep=crcr]{%
1	1\\
2	2\\
4	4\\
8	8\\
16	16\\
};
\end{axis}

\end{tikzpicture}%
    
    \caption{Strong scaling. Solve time vs number of GPUs. Time needed for the FCG to reach convergence (left panel) and obtained speedup (right panel).}
    \label{fig:strong-scaling}
\end{figure}
On a single GPU, the global memory usage—including the preconditioner and the memory required for the FCG vectors—amounts to 45 GB. This is close to the GPU's total available memory (64 GB). The results, shown in Fig.~\ref{fig:strong-scaling}, depict the solution time (left panel, solid line) compared to the ideal time (dashed black line) and the achieved speedup (right panel, solid line) relative to the optimal speedup (dashed black line). For 16 GPUs, this yields an efficiency of 44\%.

\vspace{-2em}\section{Conclusions and future perspectives}
\label{sec:conclusions}
We have demonstrated that the PSCToolkit suite effectively enables the solution of large linear
systems using a substantial number of GPUs. Notably, on the classical Poisson benchmark we show
performances comparable to, and in some configurations surpassing, the state of the art. The ease of
implementation provided by these libraries allows users to achieve high performance with minimal
effort.

Looking ahead, our focus will be on improving OpenMP and GPU supports for the preconditioners setup
and studying and implementing robust and efficient polynomial smoothers optimized for GPU use.
Additionally, we are looking at  extending support beyond NVIDIA GPUs to include those from other
manufacturers. To this end, we are monitoring the development of compilers providing OpenMP and
OpenACC support; these environments would greatly facilitate many development tasks. At the time of
writing, directive-based programming environments still exhibit a non-trivial performance penalty for
the matrix-vector product kernels of interest, whilst we have  prioritized the effective usage of our
computing platform.  To improve maintainability, we strived to confine the use of device-specific
code such as CUDA to as few lines of code as possible; an increased use of high-level compilers is the
subject of current research, even if in some cases the innermost kernels may require a  device-specific implementation.

\begin{acknowledgement}
The three authors are member of the INdAM-GNCS group. F.D. acknowledges the MIUR Excellence Department Project awarded to the Department of Mathematics, University of Pisa, CUP I57G22000700001. This work was partially supported by: Spoke 1 ``FutureHPC \& BigData'' and Spoke 6 ``Multiscale Modelling \& Engineering Applications'' of the Italian Research Center on High-Performance Computing, Big Data and Quantum Computing (ICSC) funded by MUR Missione 4 Componente 2 Investimento 1.4: Potenziamento strutture di ricerca e
creazione di ``campioni nazionali di R\&S(M4C2-19)'' - Next Generation EU (NGEU); the ``Energy Oriented Center of Excellence (EoCoE III): Fostering the European Energy Transition with Exascale" EuroHPC Project N. 101144014, funded by European Commission (EC); the ``INdAM-GNCS Project: Metodi basati su matrici e tensori strutturati per problemi di algebra lineare di grandi dimensioni'', CUP E53C22001930001; the European Union - NextGenerationEU under the National Recovery and Resilience Plan (PNRR) - Mission 4 Education and research - Component 2 From research to business - Investment 1.1 Notice Prin 2022 - DD N. 104 2/2/2022, titled ``Low-rank Structures and Numerical Methods in Matrix and Tensor Computations and their Application'', proposal code 20227PCCKZ -- CUP I53D23002280006. The experiments on the EuroHPC-JU Leonardo supercomputer have been run under the Leonardo Early Access Program, Grant LEAP\_014.
\end{acknowledgement}
\ethics{Competing Interests}{The authors have no conflicts of interest to declare that are relevant to the content of this chapter.}

\vspace{-2em}
\bibliographystyle{spmpsci}
\bibliography{bibliography.bib}
\end{document}